\documentclass[reqno,a4paper,12pt]{amsart}

\usepackage[margin= 24 mm]{geometry}

\usepackage{cite}
\usepackage{amsmath,amssymb,amsfonts,amsthm}

\usepackage{graphicx}
\usepackage{tikz}
\usepackage{varwidth}

\newlength{\mycapwidth}
\newcounter{mycapcounter}
\renewcommand\caption[1]{\refstepcounter{mycapcounter}
  \bigskip\\ \begin{varwidth}[t]{\mycapwidth}
  {\bf Figure \arabic{mycapcounter}:} \ #1 \end{varwidth}}

\newcounter{myalistcount}
\newenvironment{alist}
{  \begin{list} {$(\alph{myalistcount})$}
   {\usecounter{myalistcount}\setlength\labelwidth{10 cm}}  }
{\end{list}}
\def\bal{\begin{alist}} \def\eal{\end{alist}}

\newcounter{myilistcount}
\newenvironment{ilist}
{  \begin{list} {$(\roman{myilistcount})$}
   {\usecounter{myilistcount}\setlength\labelwidth{10 cm}}  }
{\end{list}}
\def\bil{\begin{ilist}} \def\eil{\end{ilist}}

\newenvironment{blist}
{  \begin{list} {\large$\bullet$}
   {\usecounter{myblistcount}\setlength\labelwidth{10 cm}}  }
{\end{list}}
\def\bbl{\begin{blist}} \def\ebl{\end{blist}}

\numberwithin{equation}{section}

\newcommand{\eqn}[1]{\begin{equation} #1 \end{equation}}

\newcommand{\aln}[1]{\begin{alignat}{10} #1 \end{alignat}}
\newcommand{\alns}[1]{\begin{alignat*}{10} #1 \end{alignat*}}

\let\le=\leqslant
\let\ge= \geqslant
\def \ds{\displaystyle}

\newcommand\al{\alpha}
\newcommand\be{\beta}
\newcommand\De{\Delta}
\newcommand\ep{\epsilon}
\newcommand\ga{\gamma}
\newcommand\la{\lambda}
\newcommand\si{\sigma}

\def\R{{\mathbb R}}

\newcommand\bfal{\boldsymbol{\alpha}}
\newcommand\bfbe{\boldsymbol{\beta}}
\newcommand\bfeta{\boldsymbol{\eta}}

\newcommand\bzero{{\boldsymbol{0}}}

\newcommand\D{{\mathcal D}}
\renewcommand\O{{\mathcal O}}
\renewcommand\S{{\mathcal S}}
\newcommand\bS{\overline \S}
\newcommand\C{{\mathcal C}}

\newcommand\X{\times}
\newcommand{\pa}{\partial}
\newcommand\sgn{{\rm sgn}\,}

\newtheorem{theorem}{Theorem}[section]
\newtheorem{lemma}[theorem]{Lemma}

\newtheorem{corollary}[theorem]{Corollary}

\theoremstyle{definition}
\newtheorem{definition}[theorem]{Definition}

\newtheorem{remark}[theorem]{Remark}

\renewcommand{\proof}{\medskip\noindent{\em Proof.} \ }
\def \eop {\qed\medskip}

\begin{document}

\title[Linear and nonlinear problems with multi-point boundary
conditions]{Linear and nonlinear, second-order problems with
Sturm-Liouville-type, multi-point boundary conditions}
\author{Bryan  P.  Rynne}
\address{Department of Mathematics and the Maxwell Institute for
Mathematical Sciences, Heriot-Watt University,
Edinburgh EH14 4AS, Scotland.}
\email{b.p.rynne@.hw.ac.uk}

\subjclass{34B15}
\keywords{Ordinary differential equations, nonlinear boundary value
problems, multi-point boundary conditions, linear eigenvalue problem}

\begin{abstract}
We consider the nonlinear equation
$$-u'' = f(u) + h , \quad \text{on} \quad (-1,1),$$
where $f : \R \to \R$ and $h : [-1,1] \to \R$ are continuous,
together with  general Sturm-Liouville type, multi-point boundary
conditions at $\pm 1$.
We will obtain existence of solutions of this boundary value problem
under certain
`nonresonance' conditions, 
and also Rabinowitz-type global bifurcation results, which yield 
nodal solutions of the problem.

These results rely on the spectral properties
of the eigenvalue problem consisting of the equation
$$-u'' = \lambda u, \quad \text{on} \quad (-1,1),$$
together with the multi-point boundary conditions.
In a previous paper it was shown that, under certain `optimal'
conditions, the basic spectral properties of this eigenvalue problem are
similar to those of the standard Sturm-Liouville problem with
single-point boundary conditions.
In particular, for each integer $k \ge 0$ there exists a unique, simple
eigenvalue $\lambda_k$, whose eigenfunctions have `oscillation count'
equal to $k$, where the `oscillation count' was defined in terms of a
complicated Pr\"ufer angle construction.

Unfortunately, it seems to be difficult to apply the Pr\"ufer angle
construction to the nonlinear problem. Accordingly, in this paper we use
alternative, non-optimal, oscillation counting methods to obtain the
required spectral properties of the linear problem, and these are then
applied to the nonlinear problem to yield the results mentioned above.
\end{abstract}

\maketitle

% \footnotetext{\red{}}

\section{Introduction}  \label{intro.sec}

In this paper we consider the nonlinear boundary value problem 
consisting of the equation
\begin{equation} \label{nonlinear_de.eq}
 -u'' = f(u) + h ,  \quad \text{on $(-1,1)$},
\end{equation}
where $f : \R \to \R$ and $h : [-1,1] \to \R$ are continuous,
together with the multi-point boundary conditions
\begin{equation} \label{slbc.eq}
\alpha_0^\pm u(\pm 1)  + \beta_0^\pm u'(\pm 1) =
\sum^{m^\pm}_{i=1} \alpha^\pm_i u(\eta^\pm_i)
 +
\sum_{i=1}^{m^\pm} \beta^\pm_i u'(\eta^\pm_i) ,
\end{equation}
where
$m^\pm \ge 1$ are integers,
$\alpha_0^\pm,\beta_0^\pm \in \R$,
and, for each
$i = 1,\dots,m^\pm$,
the numbers
$\alpha_i^\pm,\beta_i^\pm \in \R$,
and
$\eta_i^\pm \in [-1,1]$,
with $\eta_i^\pm \ne \pm 1$.
However, as a preliminary to this we will discuss the linear eigenvalue
problem consisting of the equation
\begin{equation} \label{eval_de.eq}
 -u'' = \la u ,  \quad \text{on $(-1,1)$},
\end{equation}
where $\lambda \in \R$,
together with the boundary conditions \eqref{slbc.eq}.
Naturally, an {\em eigenvalue} is a number $\lambda$ for which
\eqref{slbc.eq}-\eqref{eval_de.eq},
has a non-trivial solution $u$ (an {\em eigenfunction}).
The {\em spectrum}, $\sigma$,  is the set of eigenvalues.

Throughout the paper we will suppose that the coefficients in 
\eqref{slbc.eq} satisfy the conditions
\begin{gather}
\alpha_0^\pm \ge 0, \quad \alpha_0^\pm + |\beta_0^\pm| > 0 ,
\label{albe_nz.eq}
\\
\pm \beta_0^\pm \ge 0 ,
\label{albe_sign.eq}
\\
\left( \frac{\sum_{i=1}^{m^\pm} |\alpha_i^\pm|}{\alpha_0^\pm} \right)^2
 +
\left( \frac{\sum_{i=1}^{m^\pm} |\beta_i^\pm|}{\beta_0^\pm} \right)^2
 < 1 
\label{AB_lin_cond-2.eq}
\end{gather}
with the convention that if any denominator in
\eqref{AB_lin_cond-2.eq} is zero then the corresponding numerator must
also be zero, and the corresponding fraction is omitted from
\eqref{AB_lin_cond-2.eq}
(by \eqref{albe_nz.eq}, at least one denominator is
nonzero in each condition \eqref{AB_lin_cond-2.eq}).

Although the boundary conditions \eqref{slbc.eq} are non-local, for ease
of discussion we will say that the condition with superscript $\pm$ holds
`at the end point $\pm 1$', and we will denote these individual
conditions by $\eqref{slbc.eq}^\pm$
(and similarly for other conditions such as 
$\eqref{albe_nz.eq}^\pm$-$\eqref{AB_lin_cond-2.eq}^\pm$,
and others below).
Also, we will write
$\alpha^\pm := (\alpha_1^\pm,\dots,\alpha_{m^\pm}^\pm) \in \R^{m^\pm}$,
and similarly for $\beta^\pm$, $\eta^\pm$.
The notation $\alpha^\pm = 0$ or $\beta^\pm = 0,$ will mean the zero
vector in $ \R^{m^\pm }$, as appropriate.
When $\alpha^\pm = \beta^\pm = 0$ the multi-point boundary conditions
$\eqref{slbc.eq}^\pm$ reduce to standard, {\em single-point} conditions 
at $x = \pm 1$,
and the overall multi-point problem
\eqref{slbc.eq}-\eqref{eval_de.eq} reduces to
a standard, linear Sturm-Liouville problem.
Thus, we will term the conditions  \eqref{slbc.eq}
{\em Sturm-Liouville-type} boundary conditions.
If $\be_0^\pm = 0$ (respectively, $\al_0^\pm = 0$) we term the condition
$\eqref{slbc.eq}^\pm$  {\em Dirichlet-type}
(respectively, {\em Neumann-type}).
This terminology is motivated by observing that a multi-point
Dirichlet-type (respectively Neumann-type)
condition reduces to a single-point
Dirichlet (respectively Neumann)
condition when
$\alpha=0$ (respectively $\beta=0$).

Various types of boundary value problems with multi-point boundary
conditions, both linear and nonlinear, have
been extensively studied recently,
see for example,
\cite{DR,KK,KKW,KG,KKKW,MOR1,RYN3,RYN5,RYN6,RYN7,XU},
and the references therein.
In this paper we will continue the investigation of the spectral
properties of the linear problem 
\eqref{slbc.eq}-\eqref{eval_de.eq},
and then apply these properties to the nonlinear problem
\eqref{nonlinear_de.eq}-\eqref{slbc.eq}.

\subsection{Previous spectral results}

The spectral properties of the standard linear, single-point,
Sturm-Liouville problem are of course well known, 
% see for example \cite{CL},
but the spectral properties of the above general linear, multi-point
problem \eqref{slbc.eq}-\eqref{eval_de.eq}
are still being investigated.
Indeed, it is only recently that the basic spectral properties of
any multi-point problems have been obtained
(the multi-point problem is not self-adjoint, so in principle is more
`difficult' than the single-point problem, which is self-adjoint).
Initially, problems with a single-point condition at one end point and a
multi-point condition at the other end point were discussed, 
essentially using shooting from the single-point end,
see for example \cite{DR,MOR1,RYN3};
the papers \cite{KK,KKW} also discuss this case, with a variable 
coefficient function in the differential equation \eqref{eval_de.eq}.
Multi-point conditions at both end points are more difficult to deal
with.
Dirichlet-type conditions, at both end points, were discussed in
\cite{KKKW,RYN5}, while Neumann-type conditions were discussed in 
\cite{RYN6}.
The case of a  Dirichlet-type condition at one end point and a
Neumann-type condition at the other end point was also discussed in
\cite{RYN6}, where such conditions were termed {\em mixed}.
The full Sturm-Liouville type conditions \eqref{slbc.eq} were
discussed in \cite{RYN7}.

For the problems discussed in the papers \cite{DR,RYN3,RYN5,RYN6,RYN7}
it was shown that the spectra of these  problems have many of the 
`standard' properties of the spectrum of the usual single-point 
Sturm-Liouville problem,
specifically:
\begin{itemize}
\item[($\sigma$-a)]
$\sigma$ consists of a strictly increasing sequence of real eigenvalues
$\lambda_k$, $k=0,1,\dots;$
\item[($\sigma$-b)]
$\lim_{k \to \infty} \lambda_k = \infty$;
\end{itemize}
for each $k \ge 0$:
\begin{itemize}
\item[($\sigma$-c)]
$\lambda_k$ has geometric multiplicity 1;
\item[($\sigma$-d)]
the eigenfunctions corresponding to $\lambda_k$ have
an `oscillation count' equal to $k$.
\end{itemize}
In the single-point problem the oscillation count referred to in
property ($\sigma$-d) is simply the number of interior {\em nodal}
zeros of an eigenfunction in the interval $(-1,1)$.
However, in the multi-point problem it was found in \cite{RYN5} and
\cite{RYN6} that this method of counting eigenfunction oscillations no
longer necessarily yields property ($\sigma$-d), and alternative methods
were adopted, with different approaches being used for different types of
boundary conditions
(a more detailed discussion is given in  \cite[Section~9.4]{RYN6}).
The methods of \cite{RYN5,RYN6} were then unified and extended to the 
general Sturm-Liouville type boundary conditions in \cite{RYN7}, using a 
Pr\"ufer angle approach to describe the eigenfunction oscillation count, 
and the above spectral properties ($\sigma$-a)-($\sigma$-d) were 
obtained under the above hypotheses
\eqref{albe_nz.eq}-\eqref{AB_lin_cond-2.eq}.
It was also shown in \cite{RYN7} that these hypotheses 
were optimal for this result, in the sense  that if they do not hold then
properties ($\sigma$-c) or ($\sigma$-d) may not be true -- in the former
case, there may be eigenvalues with multiplicity 2, while in the latter
case, there may be `missing' eigenvalues, that is, there may be values of
$k \ge 0$ for which there is no eigenvalue whose eigenfunctions have the
corresponding oscillation count $k$.

Unfortunately, although the Pr\"ufer angle approach used in \cite{RYN7} 
works well for the linear eigenvalue problem, it does not seem to work 
so well for the nonlinear problem.
In particular, Rabinowitz-type global bifurcation results have often 
been obtained for both single-point and multi-point Sturm-Liouville
problems
(dating back to Rabinowitz' seminal paper \cite{RAB-glob_bif} for the
single-point problem and, for example, in \cite{RYN5} and \cite{RYN6}
for the Dirichlet-type and Neumann-type problems respectively), 
and these results have been used to obtain nodal solutions for these
problems.
Such results rely on the preservation of the nodal properties along the 
bifurcating continua, and this seems to be difficult to verify for the 
Pr\"ufer angle oscillation counting method used in \cite{RYN7}.
Hence, in this paper we will discuss how the oscillation counting methods
of \cite{RYN5} and \cite{RYN6} can be combined to apply to both the
linear and the nonlinear Sturm-Liouville multi-point boundary value
problems.
The results obtained in this manner will not be optimal for the linear
eigenvalue problem,
but they will yield global bifurcating continua and nodal solutions for
the nonlinear problem.

The linear results below will also clarify the use of two different
oscillation counting methods for the Dirichlet-type and Neumann-type
problems in \cite{RYN5} and \cite{RYN6}, and show that these problems can
be regarded as  extreme ends of a range of Sturm-Liouville boundary
conditions, with a gradual switch between the two oscillation
counting methods --- see Remark~\ref{switch_between_S_T.rem} for a more
careful description of this.
Of course, the Pr\"ufer angle oscillation count used in \cite{RYN7} also
subsumed and generalised the various oscillation counting methods used in
\cite{RYN5} and \cite{RYN6}  in the Dirichlet-type, Neumann-type 
(and mixed) cases respectively.

\section{Preliminaries}  \label{Preliminaries.sec}

\subsection{Some further notation}  \label{notation.subsec}

Clearly,  the eigenvalues $\lambda_k$ (and other objects to be
introduced below) depend on the values of the coefficients
$\alpha_0^\pm,\,\beta_0^\pm,\,\alpha^\pm,\,\beta^\pm,\,\eta^\pm,$
but in general we regard these coefficients as fixed,
and omit them from our notation.
However, at certain points of the discussion it will be convenient to
regard some, or all, of these coefficients as variable,
and to indicate the dependence of various functions on these
coefficients.
To do this concisely we will write:
\alns{
\bfal_0 &:= (\alpha_0^-,\alpha_0^+) \in \R^2
\quad \text{(for given numbers $\alpha_0^\pm \in \R$);}
\\
\bfal &:= (\alpha^-,\alpha^+) \in \R^{m^- + m^+}
\quad \text{(for given coefficient vectors $\alpha^\pm \in \R^{m^\pm}$);}
}
and similarly for $\bfbe_0,\,\bfbe,\,\bfeta$.
We also define $\bzero := (0,0)  \in \R^{m^- + m^+}$.
We may then write, for example, $\lambda_k(\bfal,\bfbe)$
to indicate the dependence of $\lambda_k$ on $(\bfal,\bfbe)$.

In addition, when discussing an individual boundary condition
\eqref{slbc.eq} at $x=-1$ or $x=1$,
it will be convenient to let $\nu$ denote one of the signs $\{\pm\}$,
and to use  the notation $ \eqref{slbc.eq}^\nu $ to refer
to the boundary condition \eqref{slbc.eq} at the specific end point  
$x = \nu 1$
(with the natural interpretation of this);
we will use a similar notation for other conditions, such as
\eqref{albe_nz.eq}-\eqref{AB_lin_cond-2.eq}, or other conditions below.
Also, for $u \in C^1[-1,1]$, the notation $u(\nu)$ or $u'(\nu)$  will 
denote the value of $u$ or $u'$ at the end point $x = \nu 1$.

For any integer $n \ge 0$, let $C^n[-1,1]$ denote the usual Banach
space of $n$-times continuously differentiable functions on $[-1,1]$,
with the usual sup-type norm, denoted by $|\cdot|_n$.
We now define an operator formulation of the differential
operator with multi-point boundary conditions.
Let
\begin{align*}
X &:= \{u \in C^2[-1,1] :
\text{$u$ satisfies \eqref{slbc.eq}} \},\quad \|\cdot \|_X := |\cdot |_2,
\\
Y &:=C^0[-1,1],\quad \|\cdot \|_Y := |\cdot |_0 ,
\\
\De u &:= u'',  \quad u \in X.
\end{align*}
By the definition of the spaces $X$, $Y$, the linear operator
$\De : X \to Y$ is well-defined and bounded,
and we can rewrite the eigenvalue problem 
\eqref{slbc.eq}-\eqref{eval_de.eq}
as
\begin{equation}  \label{mp_eval.eq}
 -\Delta u   = \lambda u , \quad u \in X .
\end{equation}

\subsection{Nodal sets} \label{Neu_nodal.sec}

The nodal/oscillation properties of solutions of nonlinear
Sturm-Liouville problems with single-point boundary conditions
are usually described in terms of sets of functions $u \in C^2[-1,1]$ 
having a specified number of interior zeros 
(that is, points $x \in (-1,1)$ for which $u(x)=0$),
and satisfying the given boundary conditions,
see, for example, \cite[Section 2]{RAB-glob_bif}.
However, in the case of multi-point boundary conditions it has also been 
found useful to count the interior zeros of $u'$.
Specifically, in \cite{RYN5} and \cite{RYN6} certain sets, 
denoted $T_k$ and $S_k $, 
were used to count oscillations in the Dirichlet-type and Neumann-type
cases respectively. 
We recall the definitions of these sets here.
For any $C^1$ function $u$, if $u(x_0)=0$ then $x_0$ is a
{\em simple} zero of $u$ if $u'(x_0) \ne 0$.

\begin{definition}
For any integer $k \ge 0$:
\\[ 1 ex]
$S_k^+ \subset C^2[-1,1]$ is the set of
functions
$u \in C^2[-1,1]$ satisfying the conditions:\\[1 ex]
S-(a)\ \ $u(\pm 1) \ne 0$ and $ u(-1) > 0$;\\[.5 ex]
S-(b) \ $u$ has only simple zeros in $(-1,1)$, and has exactly
$k$ such zeros.\\[1 ex]
We also define $S_k^- := -S_k^+$ and $S_k := S_k^+\cup S_k^-$.
\\[ 1 ex]
$T_k^+ \subset C^2[-1,1]$ is the set of
functions
$u \in C^2[-1,1]$ satisfying the conditions:\\[1 ex]
T-(a) \  $u'(\pm 1) \ne 0$ and $ u'(-1) > 0$;\\[.5 ex]
T-(b) \ $u'$ has only simple zeros in $(-1,1)$, and has exactly
$k$ such zeros;\\[.5 ex]
T-(c) \ $u$ has a zero strictly between each consecutive zero of 
$u'$.\\[1 ex]
We also define $T_k^- := -T_k^+$ and $T_k:=T_k^+\cup T_k^-$.
\end{definition}
\medskip

Clearly, the sets $S_k^\nu$ (respectively $T_k^\nu$),
$k \ge 0$, $\nu \in \{\pm\}$, are disjoint and open in $C^2[-1,1]$.
Another class of nodal spaces, denoted 
$P_k = P_k^+ \cup P_k^- \subset \R \X C^2[-1,1]$, $k \ge 0$, 
was defined in \cite{RYN7} to deal with the general Sturm-Liouville type
boundary conditions \eqref{slbc.eq} above.
These sets were defined in terms of the Pr\"ufer angle of solutions of 
\eqref{eval_de.eq};
the definition is quite long, and the details are not required here, so 
will be omitted.
Suffice it to say that the Pr\"ufer angle approach extends and unifies
the two separate approaches adopted in \cite{RYN5} and \cite{RYN6}.

\subsection{The nodal properties of the eigenfunctions }   
\label{spec.sec}

The eigenvalues and eigenfunctions of
\eqref{slbc.eq}-\eqref{eval_de.eq}
will be denoted by $\la_k$, $\psi_k$,  $k \ge 0$.
The eigenfunctions will always be normalised so that
$|\psi_k|_0 = 1$.
When $(\bfal,\bfbe) = (\bzero,\bzero)$ the multi-point boundary 
conditions \eqref{slbc.eq} reduce to the standard (Robin) conditions
\begin{equation}  \label{single_BC_both.eq}
\al_0^\pm u(\pm 1) + \be_0^\pm u'(\pm 1) = 0,
\end{equation}
and the eigenvalues and eigenfunctions of \eqref{eval_de.eq}, with 
these boundary conditions, will be denoted by
$ \la_k^{\bzero}, \ \psi_k^{\bzero}$,  $k \ge 0$.

The following theorem was proved in \cite[Theorem~4.8]{RYN7}.

\begin{theorem}  \label{spec.thm}
The spectrum $\sigma$ of $-\De$ consists of a strictly 
increasing sequence of real eigenvalues
$\lambda_k \ge 0$, $k=0,1,\dots,$
such that $\ds\lim_{k \to \infty} \lambda_k = \infty,$
and for each $k \ge 0\!:$
\bal
\item
$\lambda_k$ has geometric multiplicity $1;$
\item
$\lambda_k$ has an eigenfunction $\psi_k$
such that $(\lambda_k,\psi_k) \in P_k$.
\eal
In the Neumann-type case $\lambda_0 = 0$, while if
$ \alpha_0^- +  \alpha_0^+ > 0$
 then $\lambda_0 > 0$.
\end{theorem}

\noindent
{\it Proof} (sketch). \
For any $k \ge 0$ the multi-point eigenvalue
$\la_k(\bfal,\bfbe)$ and eigenfunction $\psi_k(\bfal,\bfbe)$ were 
constructed in \cite{RYN7} by continuation from the 
single-point eigenvalue $\la_k^\bzero$ and eigenfunction $\psi_k^\bzero$.
In essence, this continuation construction showed that the mappings 
\begin{equation} \label{cont_construction-maps.eq}
t \to \psi_k[t] := \psi_k(t\bfal,t\bfbe) : [0,1] \to C^2[-1,1] ,
\quad
t \to \la_k[t] := \la_k(t\bfal,t\bfbe) : [0,1] \to \R,
\end{equation}
are well-defined and continuous,
and when $t=0$ they satisfy
\begin{equation} \label{cont_construction-nodes_bounds.eq}
\psi_k[0] = \psi_k(\bzero,\bzero) = \psi_k^\bzero ,
\quad
\la_k[0] = \la_k(\bzero,\bzero) = \la_k^\bzero .
\end{equation}
The multi-point  eigenvalue and eigenfunction, $\la_k$, $\psi_k$, are
then obtained by setting  $t=1$.
The properties of $\la_k$, $\psi_k$ can be derived 
from the corresponding properties of the single-point eigenvalue and 
eigenfunction, $\la_k^\bzero $, $\psi_k^\bzero $
(which are obtained from standard Sturm-Liouville theory)
by showing that they are preserved during this continuation process, as
$t$ varies from 0 to 1.
\eop

\remark
The sign condition \eqref{albe_sign.eq} ensures that
$\la_0^\bzero > 0$
(except in the Neumann-type case, when $\la_0^\bzero = 0$),
and this positivity is preserved in the continuation.
It is shown in \cite{RYN7} that if \eqref{albe_sign.eq} does not hold
then negative eigenvalues may exist, and these may have geometric
multiplicity 2
(of course, this cannot happen in the standard, single-point problem).
It is also shown in \cite{RYN7} that if \eqref{AB_lin_cond-2.eq}
does not hold then there may be values of $k$ for which there is no
eigenvalue/eigenfunction pair $(\la_k,\psi_k) \in P_k$.
Hence, the `standard' spectral properties described in
Theorem~\ref{spec.thm} may not hold if either \eqref{albe_sign.eq} or
\eqref{AB_lin_cond-2.eq} are not satisfied.

\section{Nodal properties of eigenfunctions}

We will now ascertain the nodal properties of the multi-point 
eigenfunctions described  in Theorem~\ref{spec.thm}, in terms of the 
nodal sets $T_k$ and $S_k$
(and a further class of such sets to be introduced below),
instead of the sets $P_k$ used in \cite{RYN7}.
We begin with some preliminary results, which will form the basis of the 
discussion of nodal properties. 
We first note that it can easily be shown that for any solution 
$(\la,u)$, $\la > 0$, of \eqref{eval_de.eq} we have the elementary
`energy' equalities:
\begin{equation} \label{linear_energy.eq}
\la u(x)^2 + u'(x)^2 \equiv \la |u|_0^2 = |u'|_0^2, \quad  x \in [-1,1] .
\end{equation}

\begin{lemma}  \label{nodal_preservation.lem}
Suppose that $(\la,u)$, $\la > 0$, $u \not\equiv 0$, satisfies 
\eqref{eval_de.eq} and $\eqref{slbc.eq}^\nu$, 
for some $\nu \in \{ \pm \}$.
Then$:$
\bal
\item
if
\begin{equation} \label{nonzero_ud_condn.eq}
\al_0^\nu > 
 \sum^{m^\nu}_{i=1} |\alpha^\nu_i| 
 +
\la^{1/2}   \sum_{i=1}^{m^\nu} |\beta^\nu_i| ,
\end{equation}
then $u'(\nu) \ne 0;$
\item
if
\begin{equation} \label{nonzero_u_condn.eq}
\be_0^\nu > 
\frac{1}{\la^{1/2}}   \sum^{m^\nu}_{i=1} |\alpha^\nu_i| 
 +
 \sum_{i=1}^{m^\nu} |\beta^\nu_i| ,
\end{equation}
then $u(\nu) \ne 0$.
\eal
\end{lemma}

\proof
Suppose that $u'(\nu) = 0$.
Then it follows from $ \eqref{slbc.eq}^\nu $ and \eqref{linear_energy.eq}
that
\alns{
\al_0^\nu |u|_0  & \le  
|u|_0  \sum^{m^\nu}_{i=1} |\alpha^\nu_i| 
 +
|u'|_0  \sum_{i=1}^{m^\nu} |\beta^\nu_i| 
\le  |u|_0  \sum^{m^\nu}_{i=1} |\alpha^\nu_i| 
 +
\la^{1/2} |u|_0   \sum_{i=1}^{m^\nu} |\beta^\nu_i| ,
}
which contradicts \eqref{nonzero_ud_condn.eq}, and so proves part~$(a)$.
The proof of part~$(b)$ is similar.
\eop

\subsection{The case of one multi-point boundary condition}
\label{one_mpbc.sec}

We begin the discussion of the nodal properties of the eigenfunctions by 
considering the simpler case where we only have one
multi-point boundary condition.
Specifically, in this section we suppose that in the boundary condition 
$\eqref{slbc.eq}^-$ we have
\begin{equation}  \label{single_mp_condition.eq}
\al^- = \be^- = 0 ,
\end{equation}
that is, at $x = -1$ the boundary condition
$\eqref{slbc.eq}^-$ simply reduces to the Robin
condition $\eqref{single_BC_both.eq}^-$, 
while  we retain the multi-point boundary condition at $x=1$.
The case where we only have a  multi-point condition at $x = -1$ is 
entirely similar.

We first observe that if $u$ is a non-trivial solution of
$\eqref{single_BC_both.eq}^-$, \eqref{eval_de.eq}, 
and if $\al_0^- \be_0^- \ne 0$ then $u(-1) u'(-1) \ne 0$, 
which is consistent with $u$ belonging to either of the nodal sets $S_k$,
$T_k$, for some $k \ge 0$.
On the other hand, if 
$\al_0^- = 0$ 
(respectively $ \be_0^- = 0$)
then 
$u'(-1) = 0$
(respectively $u(-1)  = 0$),
in which case $u$ cannot belong to any set $T_k$
(respectively $S_k$), $k \ge 0$,
so one of the classes of nodal sets $S_k$, $T_k$ is of no use in this
case.
However, this problem is easily remedied by simply redefining the sets
$S_k$, $T_k$ to only include functions satisfying the boundary condition
$\eqref{single_BC_both.eq}^-$.
Such redefined sets would not be open in  $C^2[-1,1]$, but would be open
in the subset of $C^2[-1,1]$ consisting of functions satisfying
$\eqref{single_BC_both.eq}^-$,
which suffices for the arguments below (in this section).
We will not mention this special case again, but in any of the following
results in this section we will implicitly suppose that we are using the
redefined sets in this case.

Next, we introduce some further definitions.
We denote the single-point eigenvalues and eigenfunctions of
\eqref{eval_de.eq}, 
with the boundary condition $\eqref{single_BC_both.eq}^-$ at $x=-1$,
together with Dirichlet or Neumann boundary conditions at $x=1$, by: 
$ \la_k^{RD},\ \psi_k^{RD},\ \la_k^{RN},  \psi_k^{RN} , $
$k \ge 0$.
It can be verified that if $(\la,u)$, $\la > 0$, satisfies
\eqref{eval_de.eq}, $\eqref{single_BC_both.eq}^-$, 
then, for each integer $k \ge 0$,
\aln{
u \in T_{k+1} & \implies \ \la_{k}^{RN}  < \la < \la_{k+1}^{RN}  ,
\label{eval_bounds_T-sBC.eq}
\\
u \in S_k & \implies  \la_{k-1}^{RD}  < \la < \la_{k}^{RD}
\label{eval_bounds_S-sBC.eq}
}
(we define $\la_{-1}^{RD} := 0$).

% The following result simply considers individual eigenvalues.

\begin{theorem}  \label{nodal_props-single_mpbc-individual.thm}
Suppose that \eqref{single_mp_condition.eq}  holds.
Then, for any integer $k_0 \ge 0$\,$:$
\bal
\item
if $\eqref{nonzero_ud_condn.eq}^+$ holds for $\la = \la_{k_0}^{RN}$
then
\eqn{ \label{efun_nodal_props_T-single_mpbc-individual.eq}
  k \le k_0 - 1 \implies \psi_k \in T_{k+1} 
  \quad \text{and} \quad
  \la_{k}^{RN}  < \la_k  < \la_{k+1}^{RN} ;
}
\item
if $\eqref{nonzero_u_condn.eq}^+$ holds for $\la = \la_{k_0}^{RD}$
then 
\eqn{ \label{efun_nodal_props_S-single_mpbc-individual.eq}
  k \ge k_0 + 1   \implies \psi_k \in S_k
  \quad \text{and} \quad
  \la_{k-1}^{RD}  < \la_k  < \la_{k}^{RD} .
}
If $\eqref{nonzero_u_condn.eq}^+$ holds for $\la = \la_{0}^{RD}$
then \eqref{efun_nodal_props_S-single_mpbc-individual.eq} also holds for
$k = 0$, so $\psi_k \in S_k$, for all $k \ge 0$.
\eal
\end{theorem}

\proof
The proof relies on the continuation construction of the eigenvalues 
and eigenfunctions, as described in the above sketch of the proof of 
Theorem~\ref{spec.thm}.
In particular, we use the mappings \eqref{cont_construction-maps.eq}, 
with the properties \eqref{cont_construction-nodes_bounds.eq}.

We first assume that $\al_0^+ \ne 0$ and $ \be_0^+ \ne 0$.
Now, for any $k \ge 0$, the eigenfunction $\psi_k^\bzero$ satisfies the
boundary condition 
$ \eqref{single_BC_both.eq}^+ $ at $x=1$,
so by \eqref{albe_nz.eq} and \eqref{albe_sign.eq},
\begin{equation}  \label{sgn_u_ud_at_1-individual.eq}
\sgn \psi_k^\bzero(1) = - \sgn \frac{d \psi_k^\bzero}{dx}(1) \ne 0 ,
\end{equation}
from which the following additional nodal and eigenvalue interlacing
properties can be obtained,
\begin{equation}  \label{single_BC_interlacing-individual.eq}
\psi_k^\bzero \in S_k \cap T_{k+1} ,
\quad 
\la_k^{RN} < \la_k^\bzero < \la_k^{RD} < \la_{k+1}^{RN} ,
\quad k \ge 0 .
\end{equation}

\noindent
$(a)$ \
Suppose that $k \le k_0 - 1$.
It follows from \eqref{eval_bounds_T-sBC.eq}
and 
part~$(a)$ of Lemma~\ref{nodal_preservation.lem}
that, for $t \in [0,1]$,
\aln{
\psi_k[t] \in T_{k+1} & \implies 
\la_{k}^{RN}  < \la_k[t] < \la_{k+1}^{RN}  
\le \la_{k_0}^{RN} ,
\label{imp-T1-sm-individual.eq}
\\
\la_k[t] \le \la_{k_0}^{RN}
& \implies 
\psi_k'[t](1) \ne 0
\implies 
\psi_k[t]  \not \in \pa T_{k+1} .
\label{imp-T2-sm-individual.eq}
}
Also, \eqref{single_BC_interlacing-individual.eq} shows that the left
hand sides of
the implications 
\eqref{imp-T1-sm-individual.eq}-\eqref{imp-T2-sm-individual.eq} 
hold when $t=0$, so by continuity the right hand sides hold for all 
$t \in [0,1]$,
and putting $t=1$ yields
\eqref{efun_nodal_props_T-single_mpbc-individual.eq}.

\noindent
$(b)$ \
The proof of 
\eqref{efun_nodal_props_S-single_mpbc-individual.eq}
is similar, using \eqref{eval_bounds_S-sBC.eq} and 
part~$(b)$ of Lemma~\ref{nodal_preservation.lem}.
The analogues of the implications
\eqref{imp-T1-sm-individual.eq}-\eqref{imp-T2-sm-individual.eq} 
in this case are:
\aln{
\psi_k[t] \in S_{k} & \implies 
\la_{k_0}^{RD} \le \la_{k-1}^{RD}  < \la_k[t] 
< \la_{k}^{RD}   ,
\label{imp-S1-sm-individual.eq}
\\
\la_k[t] \ge \la_{k_0}^{RD} & \implies 
\psi_k[t](1) \ne 0
\implies 
\psi_k[t]  \not \in \pa S_k ,
\label{imp-S2-sm-individual.eq}
}
and \eqref{efun_nodal_props_S-single_mpbc-individual.eq} now follows from
these 
implications, as before.

The final result follows from a similar argument, but the lower bound 
$0 < \la_0[t]$ is now trivial,
so we only need to prevent $\la_0[t]$ crossing
$\la_0^{RD}$, which follows from the assumption that
$\eqref{nonzero_u_condn.eq}^+$ holds for $\la = \la_{0}^{RD}$
and part~$(b)$ of Lemma~\ref{nodal_preservation.lem}.

The cases 
$\al_0^+ = 0$,  $\be_0^+ > 0$, 
and
$\al_0^+ > 0$,  $\be_0^+ = 0$,
may be proved similarly, but a generalisation of these cases is stated,
and proved, in the following corollary, so we omit any further discussion
of these cases here.
\eop
\medskip

\begin{corollary}  \label{nodal_props-smp.cor}
$(a)$ \
If $\al_0^+ > 0$,  $\be^+ = 0 $, then 
\eqref{efun_nodal_props_T-single_mpbc-individual.eq} holds
for all $k \ge 0$.
\\
$(b)$ \
If  $\al^+ = 0 $, $\be_0^+ > 0$,  then 
\eqref{efun_nodal_props_S-single_mpbc-individual.eq} holds
for all $k \ge 0$.
\\
$(c)$ \
If $\be_0^+ \ne 0$,  
\eqref{efun_nodal_props_S-single_mpbc-individual.eq} holds
for all sufficiently large $k$.  
\end{corollary}

\proof
In case $(a)$, $\eqref{nonzero_ud_condn.eq}^+$ holds 
for all $\la > 0$, so we follow the proof of 
part~$(a)$ of Theorem~\ref{nodal_props-single_mpbc-individual.thm}.
In case $(b)$ (respectively, case $(c)$), $\eqref{nonzero_u_condn.eq}^+$
holds for all $\la > 0$  
(respectively, for sufficiently large $\la > 0$), 
so we follow the proof of part~$(b)$ of
Theorem~\ref{nodal_props-single_mpbc-individual.thm}.
\eop

Theorem~\ref{nodal_props-single_mpbc-individual.thm} deals with `most'
eigenvalues, but there can be an arbitrarily large `gap' or range of 
`intermediate' eigenvalues for which neither of the hypotheses
$\eqref{nonzero_ud_condn.eq}^+$ or $\eqref{nonzero_u_condn.eq}^+$
hold, so are not covered by this theorem.
For example, if
$$
\al_0^+ = \sqrt{2} , \quad \sum_{i=1}^{m^+} |\alpha_i^+| = 1 ,
\quad
\be_0^+ = \ep \sqrt{2.1} , \quad \sum_{i=1}^{m^+} |\be_i^+| = \ep ,
$$
then \eqref{AB_lin_cond-2.eq} holds, but
$$
\eqref{nonzero_ud_condn.eq}^+ \implies \la < \frac{\sqrt{2} -1}{\ep} 
\approx \frac{0.41}{\ep}   ,
\qquad
\eqref{nonzero_u_condn.eq}^+ \implies \la > \frac{1}{\ep(\sqrt{2.1} -1)}
\approx \frac{2.23}{\ep}   ,
$$
so it follows from the Sturm comparison theorem that if $\ep$ is
sufficiently small then an arbitrarily large number of eigenvalues 
do not satisfy either 
$\eqref{nonzero_ud_condn.eq}^+$ or $\eqref{nonzero_u_condn.eq}^+$.

We can remove this gap by strengthening condition
\eqref{AB_lin_cond-2.eq} somewhat.
Specifically, if we replace  \eqref{AB_lin_cond-2.eq} with the condition
\begin{equation} \label{AB_lin_cond-1.eq}
\frac{\sum_{i=1}^{m^\pm} |\alpha_i^\pm|}{\alpha_0^\pm} 
 +
\frac{\sum_{i=1}^{m^\pm} |\beta_i^\pm|}{\beta_0^\pm} 
 < 1 
\end{equation}
(in this section the condition $\eqref{AB_lin_cond-1.eq}^-$ holds
trivially, and is irrelevant, but will be used in the next section).
The inequalities  
$\eqref{nonzero_ud_condn.eq}^\nu$
and
$\eqref{nonzero_u_condn.eq}^\nu$
are related to the condition
$\eqref{AB_lin_cond-1.eq}^\nu$.
To clarify this relationship, let
$$
J^\pm := \left( \frac{\al_0^\pm}{\be_0^\pm} \right)^2 
% \quad  \nu \in \{ \pm \}
$$
(if, for either $\nu \in \{ \pm \}$, we have $\be_0^\nu = 0$ then we set 
$J^\nu := \infty$,
and the results below hold, with the natural interpretation of this).
We now have the following corollary of
Lemma~\ref{nodal_preservation.lem} which shows that if 
$\eqref{AB_lin_cond-1.eq}^\nu$ holds then there is no gap
between the values of $\la$ for which $\eqref{nonzero_ud_condn.eq}^\nu$
and $\eqref{nonzero_u_condn.eq}^\nu$ hold.

\begin{corollary}  \label{nodal_preservation.cor}
Suppose that $(\la,u)$, $\la > 0$, $u \not\equiv 0$, satisfies 
\eqref{eval_de.eq} and, 
for some $\nu \in \{ \pm \}$,
$\eqref{slbc.eq}^\nu$ and $\eqref{AB_lin_cond-1.eq}^\nu$ hold.
Then$:$
\bal
\item
$\la \le  J^\nu  \implies$ 
$\eqref{nonzero_ud_condn.eq}^\nu$ holds $\implies u'(\nu) \ne 0;$
\item
$\la \ge  J^\nu  \implies$ 
$\eqref{nonzero_u_condn.eq}^\nu$ holds $\implies u(\nu) \ne 0$.
\eal
\end{corollary}

Next, it is easy to verify that 
$\la_k^{RN} < \la_k^{RD} < \la_{k+1}^{RN}$, so
there exists a unique integer
$k_c \ge -1$ such that
\begin{equation}  \label{k_c_defn.eq}
\la_{k_c}^{RD} < J^+ \le \la_{k_c +1}^{RD} 
\end{equation}
(if $J^+ = 0$, we set $k_c := -1$).
These definitions and eigenvalue interlacing properties are illustrated
in Fig.~\ref{eval_interlacing.fig}.

\begin{center}
\setlength\unitlength{3 pt}
\begin{picture}(120,24)
\put(43,18){$(a) \ \ J^+$}  
\put(43,13){$\overbrace{\hphantom{\hspace{21 mm}}}^{\ }$} 
\put(68,18){$(b) \ \ J^+$}
\put(64,13){$\overbrace{\hphantom{\hspace{30 mm}}}^{\ }$} 
\put(0,8){\line(1,0){120}}
\put(12,5.5){\line(0,1){5}}
\put(27,6.9){$\bullet$}
\put(42,5.5){\line(0,1){5}}
\put(63,5.5){\line(0,1){5}}
\put(76,6.9){$\bullet$}
\put(93,5.5){\line(0,1){5}}
\put(113,5.5){\line(0,1){5}}
\put(10,0){$ \la_{k_c}^{RN} $}
\put(25,0){$  \la_{k_c}^{\bzero} $}
\put(40,0){$ \la_{k_c}^{RD}  $}
\put(60,0){$ \la_{k_c+1}^{RN}  $}
\put(75,0){$  \la_{k_c+1}^{\bzero}  $}
\put(90,0){$ \la_{k_c+1}^{RD}  $}
\put(110,0){$ \la_{k_c+2}^{RN}  $}
\end{picture}
\caption{ \label{eval_interlacing.fig}
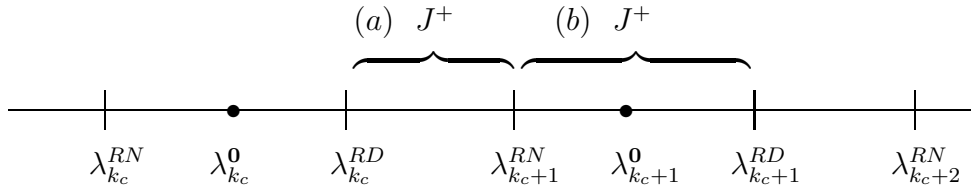
Eigenvalue interlacing and the definition of $J^+$
\\
(cases $(a)$ and $(b)$ refer to the hypotheses in
Theorem~\ref{nodal_props-single_mpbc-kc.thm}).}  
\end{center}

Combining this definition with Corollary~\ref{nodal_preservation.cor},
we see that
\alns{
\text{$\la \le \la_{k_c}^{RD}$} & \text{$\implies
\eqref{nonzero_ud_condn.eq}^+$ holds, 
}
\quad
\text{$\la \ge \la_{k_c+1}^{RD}$} & \text{$\implies
\eqref{nonzero_u_condn.eq}^+$ holds,
}
}
and combining all this with
Theorem~\ref{nodal_props-single_mpbc-individual.thm}
yields the following result.

\begin{theorem}  \label{nodal_props-single_mpbc.thm}
Suppose that \eqref{single_mp_condition.eq} and 
$\eqref{AB_lin_cond-1.eq}^+ $
hold.
Then, for any integer $k \ge 0$\,$:$
\aln{
  k \le k_c - 1 & \implies 
  \psi_k \in T_{k+1} 
  & \quad & \text{and} \quad
  \la_{k}^{RN}  < \la_k  < \la_{k+1}^{RN} ;
  \label{efun_nodal_props_T-single_mpbc.eq}
  \\
  k \ge k_c + 2  & \implies 
  \psi_k \in S_k
  && \text{and} \quad
  \la_{k-1}^{RD}  < \la_k  < \la_{k}^{RD} .
  \label{efun_nodal_props_S-single_mpbc.eq}
}
If $k_c = - 1$, that is, if $J^+ \le \la_0^{RD}$, then
\eqref{efun_nodal_props_S-single_mpbc.eq} also holds for $k = 0$, 
so $\psi_k \in S_k$, for all $k \ge 0$.
\end{theorem}

Theorem~\ref{nodal_props-single_mpbc.thm} has dealt with all the
eigenvalues in $\si$ except those with index $k = k_c$ or $k = k_c+1 $.
We deal with these in the next theorem.

\begin{theorem}  \label{nodal_props-single_mpbc-kc.thm}
Suppose that \eqref{single_mp_condition.eq} and 
$\eqref{AB_lin_cond-1.eq}^+$ hold.
Suppose, in addition, that one of the following conditions holds$:$
\bal
\item
$\la_{k_c}^{RD} < J^+ < \la_{k_c+1}^{RN}$
and either 
\smallskip
\bil
\item
$\la = \la_{k_c}^{RD}$ satisfies $\eqref{nonzero_u_condn.eq}^+$,
\smallskip
\item
$\la = \la_{k_c+1}^{RN}$ satisfies $\eqref{nonzero_ud_condn.eq}^+ ;$
\smallskip
\eil
\item
$\la_{k_c +1}^{RN} \le J^+ \le \la_{k_c +1}^{RD}. $ 
\eal
Then 
\eqref{efun_nodal_props_T-single_mpbc.eq} holds when $k = k_c$
and
\eqref{efun_nodal_props_S-single_mpbc.eq} holds when $k = k_c+1$.
\end{theorem}

\proof
The proof is similar to the proof of
Theorem~\ref{nodal_props-single_mpbc-individual.thm}.
Heuristically, we can describe the argument as follows
(again, see Fig.~\ref{eval_interlacing.fig}).
Combining the definition of $k_c$ with
Corollary~\ref{nodal_preservation.cor} shows that during the continuation
process the eigenvalues
$\la_k[t]$, $t \in [0,1]$, $k=k_c$ or $k=k_c+1$,
cannot cross either $\la_{k_c}^{RN}$ or $\la_{k_c+1}^{RD}$.
In addition, each set of hypotheses in the theorem ensures that these
eigenvalues also cannot cross one or other of 
$\la_{k_c}^{RD}$ or $\la_{k_c+1}^{RN}$.
Combining these bounds on the eigenvalues yields the result.
\eop

\remark  \label{basic_hyp_insuff-smp.rem}
By Theorem~\ref{spec.thm} above (proved in \cite[Theorem~4.8]{RYN7}),
the basic hypotheses $\eqref{AB_lin_cond-2.eq}^\pm$ 
(together with the other conditions in Section~\ref{intro.sec})
are sufficiently strong to imply that for every integer $k \ge 0$ there
is exactly one eigenvalue whose eigenfunctions lie in the nodal set $P_k$
(with either one or two multi-point boundary conditions).
It is also shown in \cite[Section~4.5]{RYN7} that if
\eqref{AB_lin_cond-2.eq} is weakened by replacing 1 on the right hand
side with $1+\ep$, for arbitrarily small $\ep > 0$, then this is no
longer true. 
However, even in the case of one multi-point condition, as considered 
in this section,  $\eqref{AB_lin_cond-2.eq}^+$ does not seem to be 
sufficiently strong to ensure that the nodal properties of all the 
eigenfunctions can always be described in terms of the nodal sets $S_k$, 
$T_k$,
for all $\al_0^+$, $\be_0^+$, $\al^+$, $\be^+$, satisfying
$\eqref{AB_lin_cond-2.eq}^+$.
On the other hand, condition $\eqref{AB_lin_cond-1.eq}^+$ does ensure 
this for all the eigenvalues, except those considered in part~$(a)$ of
Theorem~\ref{nodal_props-single_mpbc-kc.thm}.
In this case the additional conditions $(i)$ or $(ii)$ were imposed there
to deal with the eigenvalues $\la_{k_c}$ and $\la_{k_c+1}$.
These conditions represent a slight strengthening of
$\eqref{AB_lin_cond-1.eq}^+$ in that, in general, they require
$\sum^{m^+}_{i=1} |\alpha^+_i|$ 
and 
$\sum^{m^+}_{i=1} |\be^+_i|$ 
to be smaller than required by $\eqref{AB_lin_cond-1.eq}^+$.
We also note that $J^+ $ and  $k_c$ depend only on $\al_0^+$, $\be_0^+$,
so they do not `see' how small $\al^+$, $\be^+$ are, and we expect to 
obtain stronger results when these coefficient vectors are small
(when they are zero the problem reduces to the standard Sturm-Liouville
problem).

\subsection{The case of two multi-point boundary conditions}
\label{two_mpbc.sec}

We now discuss the nodal properties of the multi-point 
eigenfunctions with two multi-point boundary conditions.
We first need some more definitions.
The standard single-point eigenvalues and eigenfunctions of
\eqref{eval_de.eq}, with Dirichlet, Neumann and mixed 
(i.e., Dirichlet at one end point and Neumann at the other end)
boundary conditions, will be denoted by 
$ \la_k^D,\ \psi_k^D,\ \la_k^N,\   \psi_k^N,\ \la_k^M, \ \psi_k^M . $
It can be verified that if $(\la,u)$, $\la > 0$, satisfies 
\eqref{eval_de.eq}, then, for each integer $k \ge 0$,
\aln{
u \in T_{k+1} & \implies \ \la_k^N < \la < \la_{k+2}^N ,
\label{eval_bounds_T.eq}
\\
u \in S_k & \implies  \la_{k-2}^D < \la < \la_k^D 
\label{eval_bounds_S.eq}
}
(we define $\la_{-2}^D := 0$, $\la_{-1}^D := 0$).
We now have an analogue of 
Theorem~\ref{nodal_props-single_mpbc-individual.thm}.

\begin{theorem}  \label{nodal_props-double_mpbc-k0.thm}
For any integer $k_0 \ge 0$\,$:$
\bal
\item
if $\eqref{nonzero_ud_condn.eq}^\pm$ holds for $\la = \la_{k_0}^N$ 
then 
% \eqref{efun_nodal_props_T.eq} holds for $k;$
\eqn{  \label{efun_nodal_props_T.eq}
k \le k_0 - 2  \implies \psi_k \in T_{k+1} 
   \quad  \text{and} \quad
\la_{k}^N < \la_k < \la_{k+2}^N  ;
}
\item
if $\eqref{nonzero_u_condn.eq}^\pm$ holds for $\la = \la_{k_0}^D$ 
then 
% \eqref{efun_nodal_props_S.eq} holds for $k$.
\eqn{  \label{efun_nodal_props_S.eq}
k \ge k_0 + 2   \implies \psi_k \in S_k 
   \quad  \text{and} \quad
\la_{k-2}^D < \la_k < \la_{k}^D .
}
\eal
\end{theorem}

\proof
The proof again relies on the continuation construction of the
eigenvalues and eigenfunctions.
We first note that if 
$\al_0^\nu \ne 0$ and $\be_0^\nu \ne 0$,
for some $\nu \in \{ \pm \}$,
then, by \eqref{albe_nz.eq}, \eqref{albe_sign.eq} and
\eqref{single_BC_both.eq},
the eigenfunction $\psi_k^\bzero$ satisfies
\begin{equation}  \label{sgn_u_ud_at_nu.eq}
\sgn \psi_k^\bzero(\nu) = - \nu \, \sgn \frac{d \psi_k^\bzero}{dx}(\nu) 
\ne 0 ,
\end{equation}
from which the following nodal and eigenvalue interlacing properties can
be obtained,
\begin{equation} \label{double_BC_interlacing.eq}
\psi_k^N,\  \psi_k^{\bzero} \in S_k , \quad \psi_k^{\bzero},\ \psi_k^D 
\in T_{k+1} ,
\quad
\la_{k-1}^D = \la_k^N < \la_k^\bzero,\ \la_k^M < \la_k^D = \la_{k+1}^N ,
\quad  k \ge 0 .
\end{equation}
Also, if
$\al_0^\nu = 0$ and $\be_0^{-\nu} = 0$,
for some $\nu \in \{ \pm \}$,
then $\la_k^\bzero$ is a mixed eigenvalue, so the properties
\eqref{double_BC_interlacing.eq} again hold.
We suppose for now that either of these cases hold, and so the properties
\eqref{double_BC_interlacing.eq} hold;
the cases $\al_0^\pm  = 0$ or $\be_0^\pm = 0$, when this is not so,
will be considered below.
\medskip

\noindent
$(a)$ \
Suppose that $k \le k_0 - 2$.
Then, by \eqref{eval_bounds_T.eq} and
part~$(a)$ of Lemma~\ref{nodal_preservation.lem}, 
for $t \in (0,1]$,
\aln{
  \psi_k[t] \in T_{k+1} & \implies 
  \la_{k}^N < \la_k[t] < \la_{k+2}^N \le \la_{k_0}^N  ,
  \label{imp-T1.eq}
\\
  \la_k[t] \le \la_{k_0}^N
  & \implies 
  \psi_k'[t](\pm 1) \ne 0
  \implies 
  \psi_k[t]  \not \in \pa T_{k+1} .
  \label{imp-T2.eq}
}
Also, \eqref{double_BC_interlacing.eq} shows
that the left hand sides of the implications 
\eqref{imp-T1.eq}-\eqref{imp-T2.eq} 
hold when $t=0$, so by continuity the right hand sides hold for 
all $t \in [0,1]$,
and putting $t=1$ yields \eqref{efun_nodal_props_T.eq}.
\medskip

\noindent
$(b)$ \
If $k \ge k_0 + 2 $ then the proof of \eqref{efun_nodal_props_S.eq} is
similar.
By \eqref{eval_bounds_S.eq}, the analogues of the implications
\eqref{imp-T1.eq}-\eqref{imp-T2.eq} 
in this case are:
\aln{
  \psi_k[t] \in S_k & \implies 
  \la_{k_0}^D \le \la_{k-2}^D < \la_k[t] < \la_{k}^D ,
  \label{imp-S1.eq}
\\
  \la_k[t] \ge \la_{k_0}^D
  & \implies 
  \psi_k[t](\pm 1) \ne 0
  \implies 
  \psi_k[t]  \not \in \pa S_k ,
  \label{imp-S2.eq}
}
and \eqref{efun_nodal_props_S.eq} now follows from these implications, 
as before.

Finally, suppose that $\be_0^\pm = 0$.
Then $\la_k^\bzero$ is now a Dirichlet eigenvalue,
that is, $\la_k^\bzero = \la_k^D$,  so  although the statements in 
\eqref{double_BC_interlacing.eq} 
regarding $\la_k^\bzero$ are not all correct in this case, 
the properties of the Dirichlet and Neumann eigenvalues and 
eigenfunctions are still correct, and these suffice to show that 
the left hand sides of the implications 
\eqref{imp-T1.eq}-\eqref{imp-T2.eq} 
hold when $t=0$.
Hence, we can again obtain \eqref{efun_nodal_props_T.eq} 
by continuation.
Similarly, if $\al_0^\pm = 0$ then $\la_k^\bzero$ is a Neumann
eigenvalue  and we can again obtain \eqref{efun_nodal_props_S.eq}.
\eop

We now have the following analogue of
Corollary~\ref{nodal_props-smp.cor},
with a similar proof, based on the proof of
Theorem~\ref{nodal_props-double_mpbc-k0.thm}.

\begin{corollary}  \label{nodal_props-general-3.cor}
$(a)$ \
If $\al_0^\pm > 0$,   $\be^\pm = 0 $ then 
\eqref{efun_nodal_props_T.eq} holds for all $k \ge 0$.
\\
$(b)$ \
If  $\al^\pm = 0 $,   $\be_0^\pm > 0 $ then 
\eqref{efun_nodal_props_S.eq} holds for all $k \ge 0$.
\\
$(c)$ \
If $\be_0^\pm > 0 $  then  
\eqref{efun_nodal_props_S.eq} holds for all sufficiently large $k$.  
\end{corollary}

\remark
Corollary~\ref{nodal_props-general-3.cor} recovers the nodal
properties found in 
\cite[Theorem~5.1]{RYN5} and \cite[Theorem~5.1]{RYN6},
in the Dirichlet-type and Neumann-type cases respectively.
In fact, Corollary~\ref{nodal_props-general-3.cor} obtains slightly
more since \cite{RYN5} assumes that $\al_0^\pm = 0$, and \cite{RYN6}
assumes that $\be_0^\pm = 0$,
whereas Corollary~\ref{nodal_props-general-3.cor} allows for both
$\al_0^\pm \ne 0$ and $\be_0^\pm \ne 0$ simultaneously,
although such cases could probably have been tackled using the methods of
these previous papers.
\medskip

As in Section~\ref{one_mpbc.sec}, there is a range of intermediate 
eigenvalues not covered by 
Theorem~\ref{nodal_props-double_mpbc-k0.thm}.
We can start to deal with these using the conditions
$\eqref{AB_lin_cond-1.eq}^\pm$.
By analogy with the definition of $k_c$ in \eqref{k_c_defn.eq}, 
we define
$$
J^{\min} := \min \{ J^\pm \} , \quad  J^{\max} := \max \{ J^\pm \} ,
$$
\begin{equation}  \label{kT_kS_defn.eq}
k_T := \max \{ k : \la_k^N \le J^{\min} \} ,
\quad
k_S := \min \{ k : \la_k^D \ge J^{\max} \} .
\end{equation}
Combining these definitions with Corollary~\ref{nodal_preservation.cor}
shows that
\alns{
\text{$\la \le \la_{k_T}^{N}$} & \text{$\implies
\eqref{nonzero_ud_condn.eq}^\pm$ hold,}
\quad
\text{$\la \ge \la_{k_S}^{D}$} & \text{$\implies
\eqref{nonzero_u_condn.eq}^\pm$ hold,}
}
and combining all this with
Theorem~\ref{nodal_props-double_mpbc-k0.thm}
yields the following analogue of
Theorem~\ref{nodal_props-single_mpbc.thm}.

\begin{theorem}  \label{nodal_props-general.thm}
Suppose that $\eqref{AB_lin_cond-1.eq}^\pm $
hold.
Then, for any integer $k \ge 0\!:$
\aln{
  k \le k_T - 2  & \implies 
  \psi_k \in T_{k+1} 
  && \text{and} \quad
  \la_{k}^N < \la_k < \la_{k+2}^N  ;
%   \label{efun_nodal_props_T.eq}
\\
  k \ge k_S + 2  & \implies 
  \psi_k \in S_k 
  & \quad & \text{and} \quad
  \la_{k-2}^D < \la_k < \la_{k}^D .
%   \label{efun_nodal_props_S.eq}
}
\end{theorem}

\remark
Unfortunately, as we saw in Section~\ref{one_mpbc.sec} when dealing with
the single multi-point boundary condition case,
there is again a gap between the ranges of eigenvalues
considered in Theorems~\ref{nodal_props-double_mpbc-k0.thm}
or~\ref{nodal_props-general.thm}.
In Section~\ref{one_mpbc.sec} this gap was due to a
gap between the values of $\la$ at which the conditions 
$\eqref{nonzero_ud_condn.eq}^+$ and $\eqref{nonzero_u_condn.eq}^+$
hold, and could be eliminated by slightly strengthening
the condition $\eqref{AB_lin_cond-2.eq}^+$.
However, when we have two multi-point boundary conditions
there is also a  gap caused by the differences between the ranges of
the values of $\la$ at which the `switchover' between the conditions
$\eqref{nonzero_ud_condn.eq}^\pm$ 
and
$\eqref{nonzero_ud_condn.eq}^\pm$ 
occurs at the two end points $\pm 1$.
This gap causes a significant additional difficulty in describing the
nodal properties, which we now describe.

The proof of Theorem~\ref{nodal_props-double_mpbc-k0.thm} used the fact
that in the continuation process:
\bbl
\item
if $\la_k[t] \le \la_{k_0}^N $ then both conditions
$\eqref{nonzero_ud_condn.eq}^\pm $ hold, 
so zeros of $\psi_k'[t]$ cannot cross either of 
the end points $\pm 1$, while zeros of $\psi_k[t]$ might cross both;
\item
if $\la_k[t] \ge \la_{k_0}^D $ then  both conditions
$\eqref{nonzero_u_condn.eq}^\pm $ hold,
so zeros of $\psi_k[t]$ cannot cross either of 
the end points $\pm 1$, while zeros of $\psi_k'[t]$  might cross both. 
\ebl
In the intermediate range, 
when $\la_{k_0}^N < \la_k[t] < \la_{k_0}^D $
(which was not considered in 
Theorem~\ref{nodal_props-general.thm}),
even if $\eqref{AB_lin_cond-1.eq}^\pm $ both hold it might be the case
that 
$\eqref{nonzero_ud_condn.eq}^\nu$
holds at one end point $\nu$, while
$\eqref{nonzero_u_condn.eq}^{-\nu}$
holds at the other end point $-\nu$.
Hence, during the continuation process, a zero of $\psi_k[t]$ might cross
one end point, while a zero of $\psi_k'[t]$ might cross the other end
point, so that neither zeros of $\psi_k[t]$ nor
of $\psi_k'[t]$ are preserved during the continuation.
This renders both the classes of nodal sets $S_k$ and $T_k$, $k \ge 0$,
unsuitable for dealing with this intermediate case and necessitates the
introduction of another class of nodal sets.
We will discuss this in the following subsection.

\subsubsection{Intermediate eigenvalues when there are two multi-point
BCs}  \label{Intermediate_eigenvalues.sec}

% We now consider the nodal properties of the eigenfunctions in the 
% intermediate range not considered in 
% Theorem~\ref{nodal_props-general.thm}.
For simplicity, throughout this subsection we will suppose that
$\eqref{AB_lin_cond-1.eq}^\pm $ hold, and
$$
J^+ < J^- ;
$$
the case $ J^- < J^+ $ is similar
(the case $ J^- = J^+ $ is irrelevant in this section).
It follows from this, together with
Corollary~\ref{nodal_preservation.cor}, that if $(\la,u)$, $\la > 0$,
satisfies, $\eqref{slbc.eq}$ \eqref{eval_de.eq}, 
then
\begin{equation} \label{la_intm_u_ud_nz.eq}
J^+ \le \la \le J^- \implies u'(-1) \ne 0, \quad u(1) \ne 0 .
\end{equation} 
To utilize \eqref{la_intm_u_ud_nz.eq} we now introduce the following 
nodal sets.

\begin{definition} 
For any integer $k \ge -1$, 
$R_k^+ \subset X$ is the set of functions
$u \in X$ satisfying the following conditions:\\[1 ex]
R-(a)\ \ $u'(-1) > 0$, and
$u(1) > 0$ iff $k$ is even, $u(1) < 0$ iff $k$ is odd;
\\
R-(b) \ $u$ has only simple zeros in $(-1,1)$, and has either
$k$ or $k+1$ such zeros.
\\[1 ex]
\noindent
We also define $R_k^- := -R_k^+$ and $R_k := R_k^+\cup R_k^-$.
\end{definition} 

These sets were defined in \cite[Section 9]{RYN6},
while a motivation for their somewhat strange definition was discussed
in \cite[Section 9.4]{RYN6}.
Suffice it to say here that, combined with \eqref{la_intm_u_ud_nz.eq},
they will enable us to extend the above results to (most of) the 
intermediate eigenvalues.
It was shown in \cite[Lemma 9.2]{RYN6} that the sets $R_k^\nu$,
$k \ge -1$, $\nu \in \{\pm\}$, are disjoint and open.
In addition, if $(\la,u)$, $\la > 0$, is an arbitrary solution 
of \eqref{eval_de.eq}, then for any $k \ge 0$
\begin{equation} \label{eval_bounds_R-mBC.eq}
u \in R_k  \implies 
\la_{k-1}^M <  \la < \la_{k+1}^M  
\end{equation}
(with $\la_{-1}^M := 0$).
This is analogous to \eqref{eval_bounds_T.eq}, \eqref{eval_bounds_S.eq},
and is illustrated in Fig.~\ref{eigenfunction_limits.fig}, for the case
$k=1$.
The proof is elementary, based on the definitions, and the properties of
the sine function.
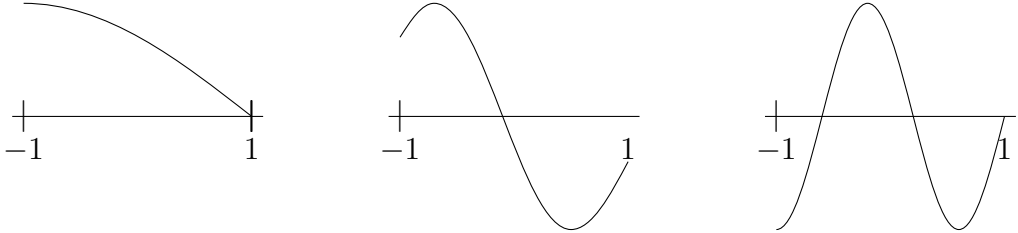
\vspace{4 mm}

\begin{center}
\newcommand \axis {\draw (\shift-1.1,0) -- (\shift+1.1,0); 
\draw (\shift+-1.0,-.1) node[below] {$-1$}; 
\draw (\shift+1.0,-.1) node[below] {$1$}; 
\draw (\shift-1.0,0) node {$|$};\draw (1.0,0) node {$|$}; 
}
\begin{tikzpicture}[scale=1.5]
\newcommand \shift {0};
\axis
\draw (-1,1) cos (1,0)  ;
\draw (0,-1.3) node[below] 
  {$\la = \la_0^M $, $u=\psi_0^M \in \pa R_1$ }; 
\renewcommand \shift {3.3};
\axis
\draw  (\shift-1,.7) sin (\shift-.7,1) cos (\shift-.1,0) sin
  (\shift+.5,-1)  cos (\shift+1,-.4)  ;
\draw (\shift+0,-1.3) node[below] 
  {$\la = \la_1^\bzero $, $u=\psi_1^\bzero \in R_1$ }; 
\renewcommand \shift {6.6};
\axis
\draw (\shift-1,-1) cos (\shift-.6,0) sin (\shift-.2,1)  cos
(\shift+.2,0) sin (\shift+.6,-1) cos (\shift+1,0) ;
\draw (\shift+0,-1.3) node[below] 
  {$\la = \la_2^M  $, $u=\psi_2^M\in \pa R_1$ }; 
\end{tikzpicture}
\vspace{-4 mm}
\caption{   \label{eigenfunction_limits.fig}
Eigenfunctions corresponding to various eigenvalues. 
}
\end{center}
\medskip

\noindent
The relevant interlacing properties are in
\eqref{double_BC_interlacing.eq}.
We also let 
\begin{equation}  \label{kT_kS_M_defn.eq}
k_{T,M} := \max \{ k : \la_k^M \le J^{\min} \} ,
\quad
k_{S,M} := \min \{ k : \la_k^M \ge J^{\max} \} .
\end{equation}
Comparing \eqref{kT_kS_M_defn.eq} with the definitions of $k_T$, $k_S$,
in \eqref{kT_kS_defn.eq}, and recalling \eqref{double_BC_interlacing.eq},
we see that 
\begin{equation}  \label{kT_kS_M_relationship.eq}
k_T-1 \le k_{T,M} \le k_T ,
\quad
k_S \le k_{S,M} \le k_S+1 .
\end{equation}

We now extend Theorem~\ref{nodal_props-general.thm} to most of the
eigenvalues omitted from that result.

\begin{theorem}  \label{nodal_props-in_gap.thm}
Suppose that $\eqref{AB_lin_cond-1.eq}^\pm $ hold.
Then, for any integer $k \ge 0,$
\begin{equation} \label{efun_nodal_props_R.eq}
k_{T,M} + 1  \le k  \le k_{S,M} - 1
\implies
\psi_k \in R_k
\quad \text{and} \quad
\la_{k-1}^M < \la_k < \la_{k+1}^M . 
\end{equation}
\end{theorem}

\proof
The proof is similar to the proof of
Theorem~\ref{nodal_props-double_mpbc-k0.thm}.
In this case, by 
\eqref{double_BC_interlacing.eq}, \eqref{la_intm_u_ud_nz.eq}
and 
\eqref{eval_bounds_R-mBC.eq},
the analogues of the implications 
\eqref{imp-T1.eq}-\eqref{imp-T2.eq} 
and
\eqref{imp-S1.eq}-\eqref{imp-S2.eq} 
are
\alns{
\psi_k[t] \in R_k & \implies 
J^+ \le \la_{k_{T,M}}^M \le \la_{k-1}^M < \la_k[t] 
  < \la_{k+1}^M \le \la_{k_{S,M}}^M \le J^- ,
% \label{imp-R1.eq}
\\
J^+ \le \la_k[t] \le J^- & \implies 
\psi_k'[t](- 1) \, \psi_k[t](1) \ne 0  
\implies 
\psi_k[t]  \not \in \pa R_k ,
% \label{imp-R2.eq}
}
and \eqref{double_BC_interlacing.eq} again shows that the left hand sides
of these implications hold when $t=0$
(a slight extension shows that $\psi_k^\bzero \in R_k$),
so \eqref{efun_nodal_props_R.eq} now follows by continuity, 
as before.
\eop
\medskip
% \vspace{4 mm}

There is still a small number (at most 4) of eigenvalues that are not
covered by
Theorems~\ref{nodal_props-general.thm} and~\ref{nodal_props-in_gap.thm},
viz., those with indices
$$
k_T - 1 \le k \le  k_{T,M}  ,
\quad  
k_{S,M} \le k \le k_S + 1 
$$
(by \eqref{kT_kS_M_relationship.eq}, each of these pairs of inequalities
corresponds to either 1 or 2 eigenvalues).
In a similar manner to the situation discussed in
Remark~\ref{basic_hyp_insuff-smp.rem}
(in the single multi-point BC case), 
the hypothesis $\eqref{AB_lin_cond-1.eq}^\pm $ does not seem to be
sufficiently strong to deal with these eigenvalues.
However, as in Theorem~\ref{nodal_props-single_mpbc-kc.thm}, a slight
strengthening of $\eqref{AB_lin_cond-1.eq}^\pm $ enables us to deal with
these eigenvalues.
We can immediately derive one such result from 
Theorem~\ref{nodal_props-double_mpbc-k0.thm}.

\begin{theorem}  \label{nodal_props-in_gap-2.thm}
Suppose that $\eqref{AB_lin_cond-1.eq}^\pm $ hold.
Then$:$
\bal
\item
if $\eqref{nonzero_ud_condn.eq}^\pm$ holds for 
$\la = \la_{k_0}^N$, with $k_0 = k_{T,M}+2$
then \eqref{efun_nodal_props_T.eq} holds for 
$k \le  k_{T,M};$
\item
if $\eqref{nonzero_u_condn.eq}^\pm$ holds for 
$\la = \la_{k_0}^D$, with $k_0 = k_{S,M}-2$
then \eqref{efun_nodal_props_S.eq} holds for 
$k \ge  k_{S,M}.$
\eal
\end{theorem}

\remark  \label{switch_between_S_T.rem}
The Dirichlet-type problem ($\be_0^\pm = 0$) was considered in
\cite{RYN5} using the sets $T_k$,
while the Neumann-type problem ($\al_0^\pm = 0$) was  considered in
\cite{RYN6} using the sets $S_k$.
Looked at in isolation this use of two diffferent types  of nodal sets
for the two cases seems slightly strange.
However, the above results now show that these cases are simply extreme
ends of a range of cases, in the following (somewhat heuristic) sense:
\bbl
\item
$\al_0^\pm = 0 \implies J^\pm = 0$: \
the nodal properties can be described using only  the sets $S_k$;
\item
$\be_0^\pm = 0 \implies J^\pm = \infty$: \
the nodal properties can be described using only  the sets $T_k$;
\item
$0 < J^\pm < \infty$: \
the nodal properties are described using a mixture of the sets
$S_k$ and $T_k$ 
(and $R_k$ in an intermediate range), 
and as $J^\pm$ varies from $0$ to $\infty$, the intermediate range
(interpreted broadly) between the sets $S_k$ and $T_k$ 
(either $k_c$ in Section~\ref{one_mpbc.sec}, or the range between
$k_T$ and $k_S$ in Section~\ref{two_mpbc.sec}) 
varies from $\infty$ to $0$.
\ebl

\section{Nonlinear problems}

\subsection{A nonresonance condition}  \label{nonres.sec}
\label{nonlin_prob.sec}

We now consider the nonlinear problem
\eqref{nonlinear_de.eq}-\eqref{slbc.eq},
which we can rewrite as
\begin{equation}  \label{nonres.eq}
 -\De u  = f(u) + h,
\quad  u \in X ,
\end{equation}
for arbitrary $h \in Y$,
where 
$f : \R \to \R$ is continuous,
and we use the notation $f : Y \to Y$ to denote the Nemitskii operator
defined by $f(u)(x) := f(u(x))$, $x \in [-1,1]$, for $u \in Y$.
We also suppose that $f$ satisfies
$\xi f(\xi) > 0$ for $\xi \in \R \setminus \{0\}$,
and
\begin{equation}  \label{f_lim_infty.eq}
0 \le f_\infty :=
\lim_{|\xi|\rightarrow\infty} \frac{f(\xi)}{\xi} \le \infty
\end{equation}
(we assume that this limit exists),

It will  be useful to  know when  $\De$ has a continuous inverse.
In the  Neumann-type case (that is, when $\alpha_0^\pm = 0$) it is clear
that any constant function $c$ lies in $X$, and $\Delta c = 0$,
so $\Delta$ cannot be invertible.
Thus, to obtain invertibility it is necessary to exclude
this case.
In view of the assumption \eqref{albe_nz.eq}
(we still assume the basic hypotheses 
\eqref{albe_nz.eq}-\eqref{AB_lin_cond-2.eq}), 
we can achieve this by imposing the further condition
\begin{equation}  \label{al_pm_strict_pos.eq}
 \alpha_0^- +  \alpha_0^+ > 0.
\end{equation}
The following results are proved in \cite[Theorem 2.1]{RYN7}
and \cite[Lemma 4.16]{RYN7}.

\begin{theorem}  \label{De_inverse.thm}
Suppose that \eqref{al_pm_strict_pos.eq} holds.
Then$:$
\bal
\item 
$\Delta : X \to Y$ has a bounded inverse $\De^{-1} : Y \to X;$
\item 
each eigenvalue $\la_k$, $k \ge 0$, is a characteristic value of the
operator $-\De^{-1} : Y \to Y$,
with algebraic multiplicity  1.
\eal
\end{theorem}

We can now obtain a solution of \eqref{nonres.eq}.

\begin{theorem}  \label{nonres.thm}
Suppose that \eqref{al_pm_strict_pos.eq} holds.
Suppose also that $f_\infty < \infty$ and $f_\infty$ is not an eigenvalue
of
$-\De$.
Then, for any $h \in Y,$ equation \eqref{nonres.eq} has a solution $u \in
X$.

\end{theorem}

\proof
The proof is similar to the proof of  \cite[Theorem~4.1]{RYN2},
using the Leray-Schauder continuation theorem,
in a relatively standard manner
(given the properties of the operator $\De^{-1}$ in
Theorem~\ref{De_inverse.thm}).
\eop

\remark
$(a)$ \
Theorem~\ref{nonres.thm} can be extended to a Sobolev space setting
(instead of the above $C^n$ setting), 
in a similar manner to that described in \cite[Remark~5.2]{DR}.
\\[1 ex]
$(b)$ \
The case where $\al_0^\pm = 0$ 
(that is, when \eqref{al_pm_strict_pos.eq} does not hold)
was considered in \cite{RYN6}, and a similar result to
Theorem~\ref{nonres.thm} was proved in \cite[Theorem~4.1]{RYN6}
(the proof is more complicated than the proof of
Theorem~\ref{nonres.thm}, due to the non-invertibility of the operator
$\De$).
Hence, we omit this case here.
\\[1 ex]
$(c)$ \
The hypothesis in Theorem~\ref{nonres.thm} that $f_\infty$ is not
an eigenvalue of $-\De$ is a `nonresonance' condition.
Nonresonance conditions have been extensively investigated
for general, single-point boundary condition problems,
both for the semilinear problem considered here, 
and the case where the linear operator $\De$ is replaced
with the quasilinear $p$-Laplacian operator $\De_p$ with $1 < p \ne 2$
(see, for example,  \cite{RYN1,RYN2} and the references therein).
For multi-point problems, 
the paper  \cite{DR} obtained similar results for a $p$-Laplacian
problem with a standard Dirichlet condition at one end-point and a
multi-point, Dirichlet-type boundary condition at the other end-point.
The results of \cite{DR} were extended to multi-point boundary conditions
at both end-points in  \cite{RYN5} (Dirichlet-type),
and in \cite{RYN6} (Neumann-type and a mixture of Dirichlet-type and
Neumann-type).
Such a nonresonance condition has not previously been obtained for the
general Sturm-Liouville-type multi-point boundary conditions
\eqref{slbc.eq}.
\medskip

\subsection{Global bifurcation theory} \label{glob_bif.sec}

In this section we consider the bifurcation problem
\begin{equation}  \label{bif.eq}
 - \De u  =\la f(u),  \quad (\la,u) \in \R \X X .
\end{equation}
where $f : \R \to \R$ is as in Section~\ref{nonres.sec},
and also satisfies
\begin{equation}  \label{f_lim_zero.eq}
0 < f_0 := \lim_{\xi \rightarrow 0} \frac{f(\xi)}{\xi} < \infty
\end{equation}
(we assume that this limit exists).
Clearly, \eqref{f_lim_zero.eq} implies that $u \equiv 0$ is a solution of
\eqref{bif.eq} for all $\la \in \R$;
such solutions will be called {\em trivial}.
We will obtain some Rabinowitz-type global bifurcation
results for the set of non-trivial solutions of \eqref{bif.eq}.

Let $\S \subset \R \X X$ denote the set of non-trivial solutions
$(\la,u)$ of \eqref{bif.eq}, and let $\bS$ denote the closure of $\S$ in
$\R \X X$.
In the following results, for any $k \ge 0$, we will use the generic
notation $N_k$ to denote one of the nodal sets
$R_k$, $S_k$ or $T_k$, and similarly for $N_k^\nu$, $\nu \in \{ \pm \}$.

\begin{lemma}\label{nhood.lem}
\bal
\item
$\bS \cap \bigl(\R\times\{0\}\bigr) \subset 
\cup_{k=0}^\infty\{(\la_k/f_0,0)\}.$
\item
Suppose that, for some $k \ge 0$, $\psi_k \in N_k$.
Then there is a neighbourhood $\O_k$ of
$(\la_k/f_0,0)$ in  $\R\X X $such that $\S \cap \O_k \subset \R \X N_k$.
\eal
\end{lemma}

\proof
Follow the proof of \cite[Lemma 4.4]{RYN1}.
\eop

For each $k \ge 0$, let  $\C_k$ denote the connected component of $\bS$
containing the point $(\la_k / f_0,0)$.
We now have the following Rabinowitz-type global bifurcation
result for the solution set of \eqref{bif.eq}.
Here, a {\em continuum} is a closed, connected set.

\begin{theorem}  \label{branches.thm}
For each $k \ge 0$ the continuum $\C_k \subset  (0,\infty) \X X  $,
and at least one of the following alternatives holds$:$
\bal
\item
$\C_k$ is unbounded in $(0,\infty) \X Y;$
\item
$(\la_j / f_0,0) \in \C_k$ for some $j \ge 0$, $j \ne k$.
\eal
\end{theorem}

\proof
Combining part~$(b)$ of Theorem~\ref{De_inverse.thm} with   
\cite[Theorem 2.3]{RAB-glob_bif} proves the result, 
with $\C_k \subset \R \X C^1[-1,1]$ in part $(a)$.
The continuity of the operator
$\De^{-1} \circ f : Y \to X$,
then shows that the continuum $\C_k$ can be regarded as a continuum in
$\R \X X$.
To show that
$\C_k \subset (0,\infty) \X X$ we note that,
by Theorem~\ref{De_inverse.thm},
the only solution $(\la,u)$ of \eqref{bif.eq} with $\la = 0$ is
$(\la,u)=(0,0)$, but by part~$(a)$ of Lemma~\ref{nhood.lem},
$(0,0) \not\in \bS$,
which implies that $\C_k \cap \big( \{0\} \X X \big) = \emptyset$.
Since $\la_k > 0$, it follows from connectedness that
$\C_k \subset (0,\infty) \X X$.
\eop
\medskip

For the problem \eqref{bif.eq} with standard
(single-point) boundary conditions, it is shown in \cite{RAB-glob_bif}
that, for each $k \ge 0$,
\begin{equation} \label{nodal_preservation.eq}
\psi_k \in N_k
\quad \text{and} \quad 
{\C_k} \backslash \{(\la_k / f_0,0)\} \subset (0,\infty) \times N_k ,
\end{equation}
with $N_k = S_k$.
That is, the nodal properties of the solutions are preserved along each
continuum $\C_k$.
Combining this with part~$(b)$  of Lemma~\ref{nhood.lem} shows that
alternative $(b)$ in Theorem~\ref{branches.thm} cannot hold for this
problem, so $\C_k$ must be unbounded 
(see \cite[Theorem 2.3]{RAB-glob_bif}).
Unfortunately, in the case of the multi-point boundary conditions
\eqref{slbc.eq}, these properties may not hold in general, but if
they do then we again obtain an unbounded continuum of solutions.
In fact, we can obtain the following result.

\begin{theorem}  \label{unbounded_branches.thm}
Suppose that, for some $k \ge 0$, 
\eqref{nodal_preservation.eq} holds for some nodal set $N_k$,
and for every $j \ge 0$, $j \ne k$,
we have $\psi_j \not\in N_k$.
Then $\C_k = \C_k^+ \cup \C_k^-$, where
\begin{equation} \label{nodal_decomposition.eq}
\C_k^\pm = 
\big( \C_k \cap ((0,\infty) \X N_k^\pm) \big)  \cup \{(\la_k/f_0,0)\} ,
\end{equation}
and each set $\C_k^\pm$ is closed, connected and unbounded in 
$(0,\infty) \X Y$.
\end{theorem}

\proof
The proof is similar to the combined proofs of
\cite[Theorem 4.5]{DR} and \cite[Theorem 4.8]{DR},
which considered the case of Dirichlet-type boundary conditions at one
end-point.
\eop

\remark
In the cases of the Dirichlet-type or Neumann-type boundary conditions
considered in \cite{DR,RYN1,RYN3,RYN5,RYN6}
it is shown that nodal properties of solutions of \eqref{bif.eq} are in
fact always preserved along the bifurcating continua, that is,
\eqref{nodal_preservation.eq} holds for all $k \ge 0$
(with $N_k = S_k$ in the Neumann-type case and $N_k = T_{k+1}$ in the
Dirichlet-type case).
Hence, for these boundary conditions
Theorem~\ref{unbounded_branches.thm} holds for all $k$,
and  we simply have the analogue of Rabinowitz' global
bifurcation theorem \cite[Theorem 2.3]{RAB-glob_bif}.
\medskip

We have not shown here how one might verify that 
\eqref{nodal_preservation.eq} holds.
We will illustrate one approach to this in the following section, and 
use the results to obtain nodal solutions of the problem.

\subsection{Nodal solutions}

We now consider the problem
\begin{equation}  \label{nodal.eq}
 -\De u = f(u), \quad  u \in X ,
\end{equation}
and we will obtain nodal solutions of this problem
(that is, solutions $u$ lying in specified nodal sets
$N_k^\nu$).
We assume that $f$ is as in Section~\ref{nonres.sec},
and the limits $f_0$ and $f_\infty$ in \eqref{f_lim_infty.eq} and 
\eqref{f_lim_zero.eq} exist.
We allow $f_\infty = \infty$, in which case we set $1/f_\infty = 0$.
Clearly, $u=0$ satisfies \eqref{nodal.eq}, so Theorem~\ref{nonres.thm} 
tells us nothing about non-trivial solutions of 
this equation.

We will use preservation of nodal properties of solutions of
\eqref{bif.eq} along the bifurcating continua, so we will need an
analogue of Lemma~\ref{nodal_preservation.lem}
for solutions of \eqref{bif.eq}.
To obtain this we first note that if $(\la,u)$ satisfies \eqref{bif.eq}
then the following generalisation of  \eqref{linear_energy.eq} can
be derived:
\begin{equation} \label{nonlinear_energy.eq}
\la F(u(x)) + u'(x)^2 \equiv \la |F(u)|_0 = |u'|_0^2 , \quad x \in
[-1,1],
\end{equation}
where $F : \R \to \R$ is defined by
$$
F(\xi) := 2 \int_0^\xi f(s)\,ds , \quad \xi \in \R .
$$
Note that for the linear equation $f(s) = s $, so that 
$F(\xi) = \xi^2$, and \eqref{nonlinear_energy.eq} reduces to
\eqref{linear_energy.eq}.
Also, by our assumption on the sign of $f$, the function $F$ is strictly 
increasing (respectively, decreasing) on 
$[0,\infty)$ (respectively, $(-\infty,0]$).
We now have the following generalisation of 
Lemma~\ref{nodal_preservation.lem}.

\begin{lemma}  \label{nonlinear_nodal_preservation.lem}
Suppose that $(\la,u)$, $\la > 0$, $u \not\equiv 0$, satisfies 
\eqref{bif.eq}.
\bal
\item
Suppose that, for some $\ga > 0$ and $\nu \in \{ \pm \}$,
\begin{equation} \label{F_small.eq}
F(\xi) \le \ga \xi^2, \quad \xi \in \R ,
\end{equation}
\begin{equation} \label{nonlinear_nonzero_ud_condn.eq}
\al_0^\nu > 
 \sum^{m^\nu}_{i=1} |\alpha^\nu_i| 
 +
\la^{1/2} \ga^{1/2}   \sum_{i=1}^{m^\nu} |\beta^\nu_i| ,
\end{equation}
then $u'(\nu) \ne 0$.
\item
Suppose that, for some $\ga > 0$ and $\nu \in \{ \pm \}$,
\begin{equation} \label{F_big.eq}
F(\xi) \ge \ga \xi^2, \quad \xi \in \R ,
\end{equation}
\begin{equation} \label{nonlinear_nonzero_u_condn.eq}
\be_0^\nu > 
\frac{1}{\la^{1/2} \ga^{1/2}}   \sum^{m^\nu}_{i=1} |\alpha^\nu_i| 
 +
 \sum_{i=1}^{m^\nu} |\beta^\nu_i| .
\end{equation}
Then $u(\nu) \ne 0$.
\eal
\end{lemma}

\proof
The proof is similar to that of Lemma~\ref{nodal_preservation.lem}, 
using $\eqref{slbc.eq}^\nu$ and
\eqref{nonlinear_energy.eq}-\eqref{nonlinear_nonzero_u_condn.eq}.
\eop

We now obtain the desired nodal solutions of \eqref{nodal.eq}.

\begin{theorem}    \label{nodal_solns_finfty_le_f0.thm}
Suppose that, for some $k \ge 0$,
\begin{equation}  \label{f_grad_infty_le_0.eq}
0 \le f_\infty < \la_k < f_0 .
\end{equation} 
\bal
\item
Suppose that $\psi_k \in T_{k+1} $, 
and $\psi_{j} \not\in T_{k+1} $ for every $j \ne k$ satisfying 
$\la_{j} \le f_0$.
Suppose also that
\eqref{F_small.eq} and $\eqref{nonlinear_nonzero_ud_condn.eq}^\pm$
hold with $\la = 1$, for some $\ga \ge f_0 $.
Then \eqref{nodal.eq} has  solutions $u_k^\pm \in T_{k+1}^\pm$.
\item
Suppose that $\psi_k \in S_k$, 
and $\psi_{j} \not\in S_k$ for every $j \ne k$ satisfying 
$\la_{j} \ge f_\infty$.
Suppose also that
\eqref{F_big.eq} and $\eqref{nonlinear_nonzero_u_condn.eq}^\pm$
hold with $\la = 1$, for some $0 < \ga \le f_\infty $.
Then \eqref{nodal.eq} has  solutions $u_k^\pm \in S_k^\pm$.
\eal
\end{theorem}

\proof
$(a)$ \ 
We note that \eqref{f_grad_infty_le_0.eq} is equivalent to
\begin{equation}  \label{equiv_crossing_condn.eq}
\la_k/f_0   < 1 < \la_k/f_\infty  .
\end{equation} 
Now, by Theorem~\ref{branches.thm} there exists a continuum $\C_k$ of
solutions of \eqref{bif.eq} bifurcating from  the point $(\la_k/f_0,0)$.
Also, by Lemma~\ref{nhood.lem} and the results of 
\cite{DAN-glob_bif} and \cite{RAB-glob_bif},
$\C_k$ can be decomposed into two subcontinua 
$\C_k = \C_k^+ \cup \C_k^-$, 
each containing $(\la_k/f_0,0)$ and such that in a neighbourhood $\O_k$ 
of $(\la_k/f_0,0)$,
\begin{equation} \label{nbhood_decomp.eq}
\big( \C_k^\pm \setminus  \{(\la_k/f_0,0)\} \big) \cap \O_k
 \cap
\big( (0,\infty) \X T_{k+1}^\pm  \big)
\ne \emptyset.
\end{equation}
To find the desired solutions of \eqref{nodal.eq} we will show that 
$\C_k^\pm$ intersects the hyperplane
$\{ 1 \} \X X$ at a pair of non-trivial solutions $(1,u_k^\pm)$ of 
\eqref{bif.eq}, 
with $u_k^\pm \in T_{k+1}^\pm$, and the functions $u_k^\pm$ are then the 
desired solutions of \eqref{nodal.eq}.
This type of argument is well-known for standard, single-point
Sturm-Liouville boundary conditions, 
see, for example, the proof of \cite[Theorem 5.3]{DR} for more
details, although the argument predates \cite{DR}.
The difficulty in the present situation is the potential non-preservation
of the nodal properties of the solutions on the continua $\C_k^\pm$.
That is, \eqref{nodal_decomposition.eq} may not hold globally,
and the continua $\C_k^\pm$ may not have the properties in 
Theorem~\ref{unbounded_branches.thm}.

To deal with this, suppose that for some $\nu \in \{ \pm \}$, 
$\C_k^\nu \cap \big( \{ 1 \} \X X \big) = \emptyset$.
% $$\C_k^\nu \cap \big( \{ 1 \} \X X_0 \big) = \emptyset ,$$
% where $X_0 := X \setminus \{0\}$.
Then, since $\la_k/f_0 < 1$ and $\C_k^\nu$ is connected, we have
$\C_k^\nu \subset (0,1) \X X$.
We will show that 
\begin{equation}  \label{C_k_in_T_k.eq}
\C_k^\nu \setminus  \{(\la_k/f_0,0)\}  \subset (0,1) \X T_{k+1}^\nu .
\end{equation} 
We first note that since any non-trivial point $(\la,u) \in \C_k^\nu$
has $\la < 1$,  the hypotheses of part~$(a)$ of
Lemma~\ref{nonlinear_nodal_preservation.lem}
hold at $(\la,u)$, which implies that $u \not\in \pa T_{k+1}^\nu$.
So, by \eqref{nbhood_decomp.eq} and the construction of the continua 
$\C_k^\pm $ in \cite{DAN-glob_bif}, if \eqref{C_k_in_T_k.eq} is false 
there must be a trivial point
$(\la_{j}/f_0 , 0) \in \C_k^\nu$, for some integer $j \ne k$ with
$\la_{j}/f_0 \le 1$ and $\psi_{j} \in T_{k+1}^\nu$.
However, this contradicts the hypothesis in the theorem, so we
conclude that \eqref{C_k_in_T_k.eq} must be true.
It now follows (similarly) from \eqref{nbhood_decomp.eq} and 
\eqref{C_k_in_T_k.eq}
that $\C_k^+ \cap \C_k^- = (\la_k/f_0,0)$, so by 
\cite[Theorem 2]{DAN-glob_bif},  $\C_k^\nu$ is unbounded.

Standard arguments (see the proof of \cite[Theorem 5.3]{DR}) now show
that there exists a sequence of non-trivial points 
$(\mu_n,v_n) \in \C_k^\nu$, $n = 1,2,\dots,$ such that,
as $n \to \infty$,
$$
\begin{cases}
\mu_n \to \la_k/f_\infty , \ | v_n |_0 \to \infty ,
& \text{if $f_\infty > 0$,}
\\
\mu_n \to \infty , 
& \text{if $f_\infty = 0$.}
\end{cases}
$$
However, each of these alternatives contradicts
\eqref{equiv_crossing_condn.eq}  and \eqref{C_k_in_T_k.eq},
so 
$\C_k^\nu \cap \big( \{ 1 \} \X X \big) \ne \emptyset$.
Next, by similar arguments to those above, it can also be shown that
there must be at least one non-trivial point in this intersection,
which completes the proof of part~$(a)$.
\\[1 ex]
$(b)$ \ In this case we use bifurcation from infinity
(see \cite{RAB-bif_inf} for more details of this) 
to obtain a continuum  $ \D_k = \D_k^+ \cup \D_k^- $ of solutions of
\eqref{bif.eq} `bifurcating from $(\la_k/f_\infty,\infty)$', 
with similar properties to those of $\C_k$ and $\C_k^\pm$.
Now, in a similar manner to the proof of part~$(a)$,
we can use  part~$(b)$ of Lemma~\ref{nonlinear_nodal_preservation.lem}
to show that if $\D_k^\nu $, $\nu \in \{ \pm \}$, does not intersect the
hyperplane $\{ 1 \} \X X$ then 
$$
(\la_k/f_0,0) \in \D_k^\nu ,
\quad
\D_k^\nu \setminus \{ (\la_k/f_0,0) \} \subset (1,\infty) \X S_k^\nu ,
$$
which again contradicts \eqref{equiv_crossing_condn.eq},
and so yields solutions  $(1,u_k^\pm) \in \{ 1 \} \X S_k^\pm$ of
\eqref{bif.eq}, and hence of \eqref{nodal.eq}.
\eop
\medskip

We can of course reverse the inequalities in 
\eqref{f_grad_infty_le_0.eq}.

\begin{theorem}    \label{nodal_solns_finfty_ge_f0.thm}
Suppose that, for some $k \ge 0$,
\begin{equation}  \label{f_grad_infty_ge_0.eq}
f_0 < \la_k < f_\infty \le \infty .
\end{equation} 
\bal
\item
Suppose that $\psi_k \in T_{k+1} $, 
and $\psi_{j} \not\in T_{k+1} $ for every $j \ne k$ satisfying 
$\la_{j} \le f_\infty$.
Suppose also that
\eqref{F_small.eq} and $\eqref{nonlinear_nonzero_ud_condn.eq}^\pm$
hold with $\la = 1$, for some $\ga \ge f_\infty $.
Then \eqref{nodal.eq} has  solutions $u_k^\pm \in T_{k+1}^\pm$.
\item
Suppose that $\psi_k \in S_k$, 
and $\psi_{j} \not\in S_k$ for every $j \ne k$ satisfying 
$\la_{j} \ge f_0 $.
Suppose also that
\eqref{F_big.eq} and $\eqref{nonlinear_nonzero_u_condn.eq}^\pm$
hold with $\la = 1$, for some $ \ga \le f_0 $.
Then \eqref{nodal.eq} has  solutions $u_k^\pm \in S_k^\pm$.
\eal
\end{theorem}

\proof
If $f_\infty < \infty$ then the proof is similar to that of 
Theorem~\ref{nodal_solns_finfty_le_f0.thm},
so will not be repeated.
If $f_\infty = \infty$ then we follow the proof of
Theorem~\ref{nodal_solns_finfty_le_f0.thm}, but modify it as in the proof
of \cite[Theorem~5.5]{DR} to obtain an unbounded sequence 
$(\mu_n,u_n) \in \C_k^\nu$, $n = 1,2,\dots,$ such that
$\mu_n \to 0 $, 
from which the result follows as before
(in this case, $\C_k^\nu$ bifurcates from $(\la_k/f_0,0)$ with
$\la_k/f_0 > 1$).
\eop

\remark
$(a)$ \
The conditions \eqref{f_grad_infty_le_0.eq}, 
\eqref{f_grad_infty_ge_0.eq}, say that the asymptotic gradients $f_0$, 
$f_\infty$ of $f$ lie on either side of the eigenvalue $\la_k$, so the 
gradient of $f$ `crosses' $\la_k$.
This type of `crossing of eigenvalues' condition is a standard condition 
used to obtain nodal solutions.
\\ $(b)$ \
In 
Theorems~\ref{nodal_solns_finfty_le_f0.thm}
and~\ref{nodal_solns_finfty_ge_f0.thm},
the conditions on the values of $\ga$ in the inequalities
\eqref{F_small.eq} and   \eqref{F_big.eq}
are related to the conditions 
\eqref{f_grad_infty_le_0.eq} 
and 
\eqref{f_grad_infty_ge_0.eq} 
on the asymptotic values of $f$.
Since
$$
\lim_{\xi \to 0} \frac{F(\xi)}{\xi^2} = f_0,
\quad
\lim_{\xi \to \infty} \frac{F(\xi)}{\xi^2} = f_\infty,
$$
we see that 
\eqref{F_small.eq} can only hold with 
$\ga \ge \max \{ f_0 , \, f_\infty \}$,
while
\eqref{F_big.eq} can only hold with 
$\ga \le \min  \{ f_0 , \, f_\infty \}$,
so the hypotheses on $\ga$, $f_0$ and $f_\infty$ in these 
theorems are consistent.
\\ $(c)$ \
Some simple sufficient conditions for
\eqref{F_small.eq} and   \eqref{F_big.eq} 
to hold, with $\ga = f_0$, are as follows:
writing $f$ in the form
$f(\xi) = \big(f_0 + g(\xi) \big) \xi , $ $\xi \in \R$,
then  
$$
\text{
$g \le 0 \implies$ \eqref{F_small.eq} holds,
\quad
$g \ge 0 \implies$ \eqref{F_big.eq} holds.
}
$$

\remark
In the proofs of 
Theorems~\ref{nodal_solns_finfty_le_f0.thm}
and~\ref{nodal_solns_finfty_ge_f0.thm}
we have obtained continua of solutions, but these do not have the full
Rabinowitz-type global properties as described in
Theorem~\ref{unbounded_branches.thm},
since the hypotheses imposed only ensure preservation of nodal properties
above or below $\la = 1$, and the continua cross $\la = 1$, so the nodal
properties need not be preserved globally.
These hypotheses could be strengthened to yield such full global results
in a variety of ways, but for brevity we will omit this here.

\end{document}